\documentclass[11pt]{amsart}

\usepackage{amsmath,amsthm,amssymb,epsfig, pstricks, pst-node,
  multirow, ifthen}

\newcounter{isamac} 
\ifthenelse{\value{isamac}=1}{
\input BoxedEPS          %
\SetTexturesEPSFSpecial  
\HideDisplacementBoxes   %
}{}

\newtheorem{theorem}{Theorem}
\newtheorem{corollary}[theorem]{Corollary}
\newtheorem{proposition}[theorem]{Proposition}
\newtheorem{lemma}[theorem]{Lemma}
\newtheorem{definition}[theorem]{Definition}
\newtheorem{remark}[theorem]{Remark}

\newtheorem{algorithm}[theorem]{Algorithm}

\newtheorem{example}[theorem]{Example}

\newcommand{\spo}{\mathrm{Sp}}      

\newcommand{\indm}[2]{\ensuremath{{\mathfrak I}_{\kern-1pt\scriptstyle#1}({\mathcal
#2})}} 

\newcommand{\ind}{\mbox{$\perp \kern-5.5pt \perp$}}
\newcommand{\nind}{\mbox{$\not\hspace{-4pt}\ind$}}

\newcommand{\uned}{\hbox{\kern3pt\raise2.5pt\vbox{\hrule
width9pt height 0.3pt}\kern3pt}}

\newcommand{\dashed}{\hbox{\kern3.05pt\raise2.5pt\vbox{\hrule
width1.7pt height 0.3pt}\kern1.8pt\raise2.5pt\vbox{\hrule
width1.7pt height 0.3pt}\kern1.8pt\raise2.5pt\vbox{\hrule
width1.7pt height 0.3pt}\kern1.8pt\raise2.5pt\vbox{\hrule
width1.7pt height 0.3pt}\kern3.05pt}}

\newcommand{\lhead}{\ensuremath{\prec}}

\newcommand{\head}{\ensuremath{\succ}}

\newcommand{\pedg}[2]{\ensuremath{{\kern0.5pt
\scriptstyle{\ifthenelse{\equal{\head}{#1}}{\lhead\kern0.5pt}{#1\kern0.5pt}}\joinrel\relbar
\negthinspace\relbar\joinrel{\kern0.5pt #2}\kern0.5pt}}}

\newcommand{\pdots}{\hbox{\kern2.5pt\raise1.5pt\hbox{\ensuremath{\ldots}}\ke
rn2.5pt}}  

\newcommand{\RRR}{\mathbb{R}}
\newcommand{\NNN}{\mathbb{N}}

\newcommand{\bi}{\leftrightarrow}

\usepackage{natbib}

\usepackage[left=3.5cm,top=3.5cm,bottom=3.5cm,right=3.5cm]{geometry}

\newcommand{\cI}{\mathcal{I}}
\newcommand{\G}{G}
\newcommand{\cC}{\mathcal{C}}
\newcommand{\bP}{\mathbf{B}}

\title[Binary Models for Marginal Independence]{Binary Models for Marginal Independence}
\author{Mathias Drton}
  \address{\em The University of Chicago, Chicago, U.S.A.}
 \author{Thomas S.\ Richardson}
  \address{\em University of Washington, Seattle, U.S.A.}
\email{{tsr@stat.washington.edu}}

\begin{document}

\begin{abstract}
  Log-linear models are a classical tool for the analysis of contingency
  tables.  In particular, the subclass of graphical log-linear models
  provides a general framework for modelling conditional independences.
  However, with the exception of special structures, marginal independence
  hypotheses cannot be accommodated by these traditional models.  Focusing
  on binary variables, we present a model class that provides a framework
  for modelling marginal independences in contingency tables.  The approach
  taken is graphical and draws on analogies to multivariate Gaussian models
  for marginal independence.  For the graphical model representation we use
  bi-directed graphs, which are in the tradition of path diagrams.  We show
  how the models can be parameterized in a simple fashion, and how maximum
  likelihood estimation can be performed using a version of the Iterated
  Conditional Fitting algorithm.  Finally we consider combining these
  models with symmetry restrictions.
\end{abstract}


\maketitle

\section{Introduction}
\label{sec:intro}

\begin{figure}[tbp]
  \begin{center}
    \small
    \newcommand{\myNode}[2]{\circlenode{#1}{\makebox[2ex]{#2}}}
    \begin{pspicture}(-1,0)(2,3)
    \psset{linewidth=0.6pt,arrowscale=1.5 2,arrowinset=0.1}  
    \psset{fillcolor=lightgray, fillstyle=solid}          
    \rput(-1,-.1){(a)}
    \rput(0, 2){\myNode{1}{$A_1$}}
    \rput(0, 0){\myNode{2}{$D_1$}} 
    \rput(2, 2){\myNode{3}{$A_2$}} 
    \rput(2, 0){\myNode{4}{$D_2$}} 
     \psset{fillcolor=white, fillstyle=solid}    
    \rput(-0.5, 1){\myNode{5}{$S_1$}}
    \rput(2.5, 1){\myNode{6}{$S_2$}} 
    \rput(1,0.5){\myNode{8}{$G_D$}} 
    \rput(1, 1.5){\myNode{7}{$G_A$}} 
    \ncline{->}{5}{1}
     \ncline{->}{5}{2}
      \ncline{->}{6}{3}
     \ncline{->}{6}{4}
     \ncline{->}{7}{1}
     \ncline{->}{7}{3}
      \ncline{->}{8}{2}
     \ncline{->}{8}{4}
    \end{pspicture}
\hskip45pt
    \begin{pspicture}(-1,0)(2,3)
    \psset{linewidth=0.6pt,arrowscale=1.5 2,arrowinset=0.1}  
    \psset{fillcolor=lightgray, fillstyle=solid}          
    \rput(-1,-.1){(b)}
    \rput(0, 2){\myNode{1}{$A_1$}}
    \rput(0, 0){\myNode{2}{$A_2$}} 
    \rput(2, 2){\myNode{3}{$D_1$}} 
    \rput(2, 0){\myNode{4}{$D_2$}} 
     \psset{fillcolor=white, fillstyle=solid}
     \rput(1, 1){\myNode{5}{$G$}} 
     \rput(-0.5, 1){\myNode{6}{$S_1$}}
    \rput(2.5, 1){\myNode{7}{$S_2$}} 
    \ncline{->}{5}{1}
    \ncline{->}{5}{2}
    \ncline{->}{5}{3}
    \ncline{->}{5}{4}
    \ncline{->}{6}{1}
     \ncline{->}{6}{2}
      \ncline{->}{7}{3}
     \ncline{->}{7}{4}
    \end{pspicture}
  \end{center}
\vskip5pt
  \caption{Possible generating models. Observed variables are shaded.
  (a) Separate genes relating to Alcohol ($G_A$) and Depression ($G_D$); (b)
  a common gene ($G$). $S_j$ represents the personal experiences of twin $j$.
  Unobserved variables are hypothesized to be independent.}
  \label{fig:generating}
\end{figure}
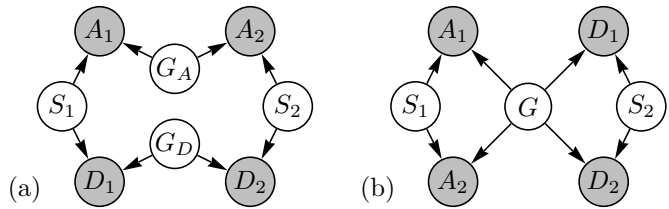

In seminal work \citet{anderson:1969,anderson:1970,anderson:1973} studied
Gaussian models defined by hypotheses that are linear in covariances.  Such
hypotheses include as a special case, zero restrictions on covariance
matrices.  These restrictions correspond to marginal independences, which
may arise for example through confounding effects of unobserved variables
\citep{coxwerm:lindep,coxwerm:book,pearl:1994,richardson:2002}.  For a
graphical representation of zero restrictions on covariance matrices,
\cite{coxwerm:lindep,coxwerm:book} introduced {\em covariance graphs\/}:
each variable is represented by a vertex; two vertices are linked by a
dashed edge
if the model does not set the corresponding covariance to zero.  Dashed
edges differentiate these graphs from undirected graphs, which represent
zero hypotheses on the inverse covariance matrix \citep{lau:bk}.  More
recently, a number of authors have used bi-directed edges
($\leftrightarrow$) in place of dashed edges which is consistent with
Sewall Wright's (\citeyear{wright:cnc}) path diagram notation; compare
Figures \ref{fig:generating}, \ref{fig:bg} and \ref{fig:trustgraph}(a)
below.  Covariance graph models have appeared in several different contexts
\citep[e.g.,][]{butte:2000,diaconis:02,grzebyk:2004,mao:2004}.  Maximum
likelihood (ML) estimation and likelihood ratio (LR) tests in these
Gaussian models can be carried out using the Iterative Conditional Fitting
algorithm \citep{drton:2003b,drton:2004jrssb}, which is implemented in the
`ggm' package in R \citep{marchetti:2006}.

There have been several efforts aimed at developing binary models with
analogous independence structure.  Kauermann (\citeyear{kauermann:1997})
uses the multivariate logistic (m-logit) transformation due to
\cite{mccullagh:89,mccullagh:89b,glonek:1995}, which consists of selecting
the highest order interaction term from every margin \citep[see
also][]{bergsma:2002}.  Cox's (\citeyear{cox:isi}) assumes that the joint
distribution is quadratic exponential, and then approximates marginal
distributions via series expansions.  An alternative approach is to use the
nonparametric concept of independence in order to form models for
categorical data that are analogous to Gaussian models.  Many existing
discrete models, such as the popular graphical log-linear models for
modelling conditional independence in contingency tables are often
motivated this way \citep{wermuth:76,darroch:80}.  In this paper we take
this route to developing
a general framework for modelling marginal independence that is a natural
counterpart to graphical log-linear models.

For an example of a marginal independence pattern that cannot be
represented using log-linear models but that our new models can accomodate
very naturally suppose that we are investigating the relationship between
alcohol dependence and depression. We have data from female mono-zygotic
twins, indicating whether or not each twin is alcohol dependent ($A_i$) and
whether or not they suffer from major depression ($D_i$); see Table
\ref{tab:mzdata}.  Consider the two graphs shown in Figure
\ref{fig:generating}.  Both hypothesize that for each twin there are
independent factors relating to individual experiences ($S_i$) which
influence both alcoholism and depression; however, graph (b) hypothesizes
in addition that there is a single genetic factor which influences both
traits, while graph (a) supposes that there is no such single factor, and
that $G_A$, $G_D$, $S_1$ and $S_2$ are mutually independent.  Graph (b)
does not imply any independence restrictions relating the observed
variables, while graph (a) implies that
\begin{equation}
  \label{eq:jointind}
  A_1 \ind D_2 \quad\hbox{ and }\quad A_2 \ind D_1;
\end{equation}
using the notation of \cite{dawid:condind}. Under (a) one twin's alcohol
dependence status is independent of the other twin's depression status.
Note that we do not make any assumption concerning the marginal
distributions of the unobserved variables.  In particular, testing the
hypothesis (\ref{eq:jointind}) provides a way of testing the scientific
hypothesis leading to graph (a) without having to specify the number of
levels of the possibly complex genetic factors.  This focus on implied
independences is in the spirit of the work on ancestral graphs
\citep{richardson:2002} and summary graphs \citep{coxwerm:book}.  We remark
that \cite{ekholm:tr:06b,ekholm:tr:06a} recently fit latent class models to
twin data including those in Table \ref{tab:mzdata}.  The precise
relationship between latent class models and the marginal independence
models we discuss in the sequel is an open problem, but if two such models
can be shown to coincide then the EM algorithm provides an alternative
method for model fitting.  However, in this context it should be noted that
there exist Gaussian covariance graph models that cannot be parameterized 
by latent variable models \citep[][\S8.6]{richardson:2002}.

\begin{table}
   \caption{Data on $n=597$  pairs of twins; adapted from \cite{kendler:92}.
  \label{tab:mzdata}}
  \centering
  \begin{tabular}{*{7}{c}}
    \hline
    & & \multicolumn{2}{c}{$D_1=0$} & & \multicolumn{2}{c}{$D_1=1$}\\
    \cline{3-4} \cline{6-7}
    & & \multicolumn{1}{c}{$D_2=0$} & \multicolumn{1}{c}{$D_2=1$} & &
    \multicolumn{1}{c}{$D_2=0$} & \multicolumn{1}{c}{$D_2=1$}\\ 
    \hline
    \multirow{2}{*}{$A_1=0$} & $A_2=0$ & 288 & 80 && 92 & 51\\
    & $A_2=1$ & 15  &  9 &&  7 & 10\\
    \hline
    \multirow{2}{*}{$A_1=1$} & $A_2=0$ &  8  &  4 &&  8 &  9\\
    & $A_2=1$ &  3  &  2 &&  4 &  7\\
    \hline
  \end{tabular}
\end{table}

If the variables were jointly Gaussian, then hypothesis (\ref{eq:jointind})
would restrict the appropriate two entries in the covariance matrix to
zero.  Hence, a likelihood ratio test of (\ref{eq:jointind}) could be
performed by fitting the covariance matrix subject to this restriction.
However, when the variables are binary, performing such a test is not at
all straightforward.  In particular, there does not exist a log-linear
model that is equal to the family of binary distributions obeying
(\ref{eq:jointind}).  In fact, the marginal independence restrictions
(\ref{eq:jointind}) correspond to complicated non-linear restrictions on
the parameters of the log-linear expansion of the joint density of $(A_1,
A_2, D_1, D_2)$.  The difficulty encountered here is an instance of the
problem of lack of compatibility of margins in log-linear parametrizations
\citep[p.534]{glonek:1995}; see also
\cite{mccullagh:89}.  In this simple example, a practical solution might be
to combine separate marginal tests, but there would be an obvious loss of
efficiency in so doing.  The methods developed in this paper allow the loss
of efficiency to be avoided by providing models that capture precisely
hypotheses like (\ref{eq:jointind}).  The fitting algorithm we present
allows tests that make use of all data available, such as LR- and
$\chi^2$-tests, to be performed.

The remainder of the paper is organized as follows.  In 
\S\ref{sec:bidi} we describe the graphical representation of marginal
independence patterns.  This representation facilitates the understanding
of marginal independence structures in multivariate normal distributions
and provides the basis for our transfer of model structure to the binary
case.  This transfer yields models that are defined implicitly in terms of
independence constraints.  In \S\ref{sec:binary-marg-ind-models}, we show
that a linear change of coordinates leads to a surprisingly simple
characterization of marginal independence.  This characterization
immediately yields a multilinear model parameterization.  ML estimation in
the proposed models is discussed in \S\ref{sec:mle} and the Iterative
Conditional Fitting algorithm for computing ML estimates is developed in
\S\ref{sec:icf}.  In \S\ref{sec:survey}, the methodology is illustrated in
an application to survey data.  In the twin data example mentioned above,
symmetry under permuting the labels 1 and 2 given to the twins is an
interesting hypothesis.  Combining such symmetry constraints with marginal
independence is the topic of section \S\ref{sec:symm-indep}.  We conclude
in \S\ref{sec:conclu}, where connections to other work are discussed.

\section{Bi-directed graphs and marginal independence}
\label{sec:bidi}

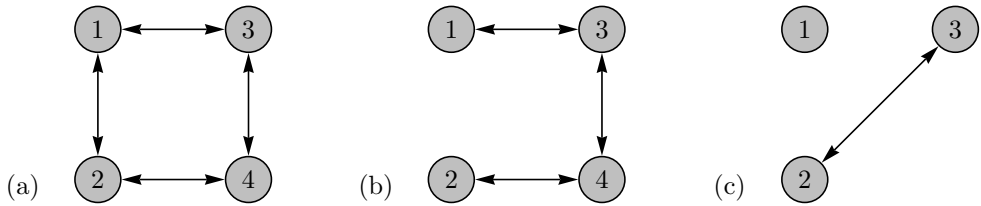
\begin{figure}[tbp]
  \begin{center}
    \small
    \newcommand{\myNode}[2]{\circlenode{#1}{\makebox[2ex]{#2}}}
    \begin{pspicture}(-1,0)(2,3)
    \psset{linewidth=0.6pt,arrowscale=1.5 2,arrowinset=0.1}  
    \psset{fillcolor=lightgray, fillstyle=solid}          
    \rput(-1,-.1){(a)}
    \rput(0, 2){\myNode{1}{$1$}}
    \rput(0, 0){\myNode{2}{$2$}} 
    \rput(2, 2){\myNode{3}{$3$}} 
    \rput(2, 0){\myNode{4}{$4$}} 
    \ncline{<->}{1}{3}
    \ncline{<->}{3}{4}
    \ncline{<->}{4}{2}
    \ncline{<->}{2}{1}
    \end{pspicture}
\hskip45pt
    \begin{pspicture}(-1,0)(2,3)
    \psset{linewidth=0.6pt,arrowscale=1.5 2,arrowinset=0.1}  
    \psset{fillcolor=lightgray, fillstyle=solid}          
    \rput(-1,-.1){(b)}
    \rput(0, 2){\myNode{1}{$1$}}
    \rput(0, 0){\myNode{2}{$2$}} 
    \rput(2, 2){\myNode{3}{$3$}} 
    \rput(2, 0){\myNode{4}{$4$}} 
    \ncline{<->}{1}{3}
    \ncline{<->}{3}{4}
    \ncline{<->}{4}{2}
    \end{pspicture}
\hskip45pt
    \begin{pspicture}(-1,0)(2,3)
    \psset{linewidth=0.6pt,arrowscale=1.5 2,arrowinset=0.1}  
    \psset{fillcolor=lightgray, fillstyle=solid}    
    \rput(-1,-.1){(c)}
    \rput(0, 2){\myNode{1}{$1$}}
    \rput(0, 0){\myNode{2}{$2$}} 
    \rput(2, 2){\myNode{3}{$3$}} 
    \ncline{<->}{2}{3}
    \end{pspicture}
  \end{center}
\vskip5pt
  \caption{(a) a bi-directed four cycle; (b) a bi-directed four chain; (c)
    graph with two disconnected components.}
  \label{fig:bg}
\end{figure}

A bi-directed graph $\G=(V,E)$ is a graph whose edges satisfy\ $(v,w)\in E$
if and only if $(w,v)\in E$. The edges are drawn bi-directed as $v\bi w$ if
$(v,w)\in E$, see Figure \ref{fig:bg}.  Bi-directed graphs are special
cases of the ancestral graphs considered in \cite{richardson:2002} and the
acyclic directed mixed graphs studied in \cite{richardson:2003}; see also
\citet[p.146]{pearl:2000}.  If a vertex $w$ is equal or adjacent to another
vertex $v$ in a bi-directed graph, then $w$ is said to be a {\em spouse\/}
of $v$, and we write $w\in\spo(v)$. For a set $A\subseteq V$, we define
$\spo(A)=\cup (\spo(v) \mid v\in A)$.  Note that $A\subseteq\spo(A)$ under
this convention.

In graphical modelling, the Markov properties of a graph, i.e., 
independence statements associated with the graph, are used to define
independence models for a random vector $X=(X_v\mid v\in V)$ whose index
set is identified with the vertex set $V$ of the graph.  The independence
models associated with bi-directed graphs are based on marginal
independence, which is manifested in the {\em connected set Markov
  property\/} of \citet[ \S4]{richardson:2003}.  A vertex set $C\subseteq V$
is {\em connected\/} if every pair of vertices $v,w\in C$ are joined by a path
on which every vertex is in $C$.  The distribution of a random vector
$X=(X_v\mid v\in V)$ is said to satisfy the connected set Markov property
if
\begin{equation}
  \label{eq:connectedMP}
  X_C \ind X_{V\setminus \spo(C)},
\end{equation}
whenever $\emptyset \neq C\subseteq V$ is a connected set.  Algorithm E in
\citet[p.~354]{knuth:68} computes equivalence classes from a list of known
equivalent pairs. This can be used to find the inclusion maximal connected
sets in a given graph by letting the edges in the graph define the
equivalent pairs.


A more exhaustive Markov property is the global Markov property, which
requires all the marginal independences in (\ref{eq:connectedMP}), but also
additional conditional independences.  More precisely, the distribution of
$X$ satisfies the global Markov property of $\G$ if
\begin{equation}\label{eq:globalMP}
A \hbox{ is separated from } B \hbox{ by } V\setminus
(A\cup B\cup C)\hbox{ in } G \hbox{ implies }X_A\ind X_B\mid X_C.
\end{equation}
Here, $A$, $B$ and $C$ are disjoint subsets of $V$, and $C$ may be empty.
The separation in (\ref{eq:globalMP}) is the usual graph-theoretic
separation in which two sets $A,B\subset V$ are separated by a third set
$D\subset V$ if any path from a vertex in $A$ to a vertex in $B$ contains a
vertex in $D$.  Despite the global Markov property being more exhaustive, a
distribution satisfies the global Markov property if and only if it
satisfies the connected set Markov property.  Completeness of the global
Markov property for bi-directed graphs follows from the completeness
results for ancestral graphs \cite[Thm.\ 7.6]{richardson:2002}.  Note also
the duality between (\ref{eq:globalMP}) and the global Markov property for
undirected graphs \cite[p.\ 32]{lau:bk}.

\begin{example}\rm ({\it Four-cycle}).
  The bi-directed graph depicted in Figure \ref{fig:bg}(a) represents the
  two pairwise independence relations:
  \[
  X_1 \ind X_4 \hbox{ and } X_2 \ind X_3
  \]
  under both the connected set, and global Markov properties.  This graph
  represents the independence hypothesis considered in the introductory
  example; compare (\ref{eq:jointind}).
\end{example}

\begin{example}\rm({\it Bi-directed four-chain}).
  Consider the bi-directed graph depicted in Figure \ref{fig:bg}(b). The
  connected set Markov property states
 \[
    X_1\ind (X_2,X_4), \quad\quad X_2\ind (X_1,X_3), \quad\quad X_3 \ind X_2, \quad\quad X_4\ind
    X_1. 
  \]
  The global Markov also states, for example, $X_1\ind X_2$,
  $X_1\ind X_2\mid X_4$, and $X_2 \ind X_3\mid X_1$.
\end{example}

Clearly, every singleton $\{v\}$ is a connected set and thus
(\ref{eq:connectedMP}) requires that 
\begin{equation}\label{eq:local}
X_v\ind X_{V\setminus \spo(v)}.
\end{equation}
It follows that if the distribution of $X$ satisfies the connected set
Markov property, then it satisfies the pairwise Markov property which
requires that $X_v\ind X_w$ whenever $v\not\bi w$. The converse is true for
multivariate normal distributions \citep[Prop.\ 2.2]{kauermann:dual} but 
false in general. We note that the Markov property in
(\ref{eq:local}) occurs in the combinatorial result known as the Lov\'{a}sz
Local Lemma \citep{erdos:lovasz:local:1975}.  In the next section we define
models using the connected set (or equivalently the global) Markov
property, and not the much less restrictive pairwise Markov property
\citep[see also][]{haber:pairwise:1986}.

\begin{example}\rm ({\it Graph with two disconnected components}).
  The pairwise Markov property for the graph in Figure \ref{fig:bg}(c)
  requires $X_1\ind X_2$ and $X_1\ind X_3$, whereas the global and
  connected set Markov property also require the stronger condition that
  $X_1\ind (X_2,X_3)$.  For example, consider the distribution of
  $(X_1,X_2,X_3)$ given by
  \begin{eqnarray*}
    &p_{000} = 0.02, \quad   p_{010} = 0.03,\quad p_{100} = 0.05, \quad p_{110} = 0.10,&\\
    &p_{001} = 0.08, \quad p_{011} = 0.12, \quad p_{101} = 0.25, \quad p_{111} =
    0.35,&
  \end{eqnarray*}
  where $p_{i_1i_2i_3}=P(X_1=i_1,X_2=i_2,X_3=i_3)$.  Then
  $X_1\ind X_2$, $X_1\ind X_3$ but $X_1\nind (X_2,X_3)$.
\end{example}

We conclude this discussion of Markov properties with a lemma that provides
a useful characterization of joint distributions of discrete random vectors
that obey the connected set Markov property.  The lemma is based on the
fact that every set $D\subseteq V$ that is not connected in $\G$ can be
partitioned uniquely into inclusion-maximal connected sets
$C_1,\ldots,C_r$,
\begin{equation}
  \label{eq:partition}
  D=C_1\dot\cup C_2\dot\cup \cdots \dot\cup C_r.
\end{equation}
Here, the symbol $\dot\cup$ denotes a union of disjoint sets.

\begin{lemma}\label{lem:conncheck}
  Let $X=(X_v\mid v\in V)$ be a discrete random vector $X=(X_v\mid v\in V)$
  taking values in the set $\cI$. The joint distribution of $X$ satisfies
  the connected set Markov property for a bi-directed graph $\G=(V,E)$ if
  and only if for every disconnected set $D\subseteq V$ it holds that
  \begin{equation}
    \label{eq:product}
    P(X_D=i_D)=\\
    P(X_{C_1}=i_{C_1})P(X_{C_2}=i_{C_2})\cdots
    P(X_{C_r}=i_{C_r}), \qquad  \forall i\in \cI,
  \end{equation}
  where $C_1,\ldots,C_r$ are the inclusion-maximal connected sets
  satisfying (\ref{eq:partition}).
\end{lemma}
\begin{proof}
  If $P$ satisfies the connected set Markov property, then it also
  satisfies the global Markov property, from which we can deduce complete
  independence of the subvectors associated with the $r$ connected sets in
  (\ref{eq:partition}),
  \begin{equation}\label{eq:disconnectind}
    X_{C_1} \ind X_{C_2} \ind \ldots \ind X_{C_r}.
  \end{equation}
  This complete independence clearly implies (\ref{eq:product}).
  
  Conversely, let $C$ be a connected set. Then $D=C\dot\cup (V\setminus
  \spo(C))$ is a disconnected set, and (\ref{eq:product}) implies in
  particular $X_C\ind X_{V\setminus
  \spo(C)}$, which is (\ref{eq:connectedMP}).
\end{proof}

\section{Binary marginal independence models}
\label{sec:binary-marg-ind-models}

Let $X=(X_v\mid v\in V)$ be a random vector with binary components, i.e., 
$X$ takes on values in the set $\cI=\{0,1\}^V$, and let $P$ be the joint
distribution of $X$. (Note that to keep notation simple, we will often use
the same letter to indicate both a set and its cardinality.)  For
$i=(i_v\mid v\in V)\in \cI$, let
\begin{equation}\label{eq:jcell}
  p_i=P(X_v=i_v  \mbox{ for all } v\in V)
\end{equation}
be the {\em joint cell probability of $i$\/}.  The multivariate Bernoulli distribution of $X$ is determined by the vector
\begin{equation}\label{eq:pvec}
  p=\big(p_i\mid i \in \cI\big)
\end{equation}
in the $2^V-1$ dimensional probability simplex $\Delta$.

Using the Markov properties discussed in the previous section we can
associate an independence model with a bi-directed graph $\G=(V,E)$.
\begin{definition}\label{def:bibimodel}
  The binary bi-directed graph model associated with $\G$ is  
  defined as the family $\bP(\G)$ of probability
  distributions for a binary random vector $X=(X_v\mid v\in V)$ that obey the
  connected set Markov property (\ref{eq:connectedMP}) for $\G$.   
\end{definition} 

We begin our study of the implicitly defined model $\bP(\G)$ by making a
change of coordinates in the probability simplex.  For $\emptyset \neq
A\subseteq V$, we call
\begin{equation}\label{eq:jzero}
q_A=P(X_A=0)=P(X_v=0  \mbox{ for all } v\in A)
\end{equation}
the {\em M\"{o}bius parameter associated with $A$\/}.  If desired, $q_A$
can be viewed as a moment for indicator variables associated with the
designated levels of the considered binary variables, namely,
\[
q_A= E\big(\textstyle\prod_{i\in A} 1_{\{X_i=0\}}\big).
\]
The $2^V-1$ M\"obius
parameters can be computed from the joint cell probabilities $p$ by the
obvious summations
\begin{equation}
  \label{eq:sum}
  q_A=\sum_{i\in\cI: i_A=0} p_i,
\end{equation}
where $i_A=(i_v\mid v\in A)$.  The summations (\ref{eq:sum}) define a map
$\mu:\Delta \to \RRR^{2^V-1}$ taking the vector of joint cell probabilities
  $p\in \Delta$ to the vector of M\"obius parameters 
\begin{equation}
  \label{eq:q}
  q=\big(q_A\mid \emptyset \neq A\subseteq V\big).
\end{equation}
We call the image $Q=\mu(\Delta)$ the {\em M\"{o}bius simplex\/}.  This
simplex has the $2^V$ vertices $t^{(A)}$, $A\subseteq V$, where for
$\emptyset\not= B\subseteq V$ the $B$-th component of $t^{(A)}$ is equal to
\[
t^{(A)}_B=
\left\{
\begin{array}{cc}
  1 & : \quad B\subseteq A,\\
  0 & : \quad B\not\subseteq A.
\end{array}
\right.
\]
Clearly, $t^{(A)}$ is the image under $\mu$ of the distribution placing
point mass on the cell $(0_A, 1_{V\setminus A})$.

\begin{proposition}\label{prop:satpara}
  The linear map 
  \begin{equation}
    \label{eq:satmap}
    \begin{array}{cccl}      
      \mu :& \Delta& \to& Q\\
      &p& \mapsto& q=(q_A\mid \emptyset \neq A\subseteq V)
    \end{array}
  \end{equation}
  is bijective. Its inverse $\nu=\mu^{-1}: Q\to \Delta$ recovers the
  joint cell probabilities as alternating sums of M\"{o}bius parameters.
  Setting $q_\emptyset=1$ we have
  \[
  p_{0_A1_{V\setminus A}}=P(X_A = 0, X_{V\setminus A} = 1) =
  \sum_{B:A\subseteq B} (-1)^{|B\setminus A|} q_B.
  \]
\end{proposition}
\begin{proof}
  By definition of $Q$, the map $\mu$ is surjective.  In order to verify
  injectivity and the claimed form of the inverse, define two functions
  $\Phi$ and $\Psi$ on the power set of $V$.  Let $\Phi (A) = q_{V\setminus
    A}$ for $A \subset V$ and $\Phi(V)=1$, and $\Psi (A) = P(X_{V\setminus
    A} = 0,X_A = 1)$.  Then $\Phi(A) = \sum_{B : B\subseteq A} \Psi(B)$ and
  the claim follows from the M\"obius Inversion Lemma \cite[p.239]{lau:bk}.
\end{proof}

The maps $\nu$ and $\mu$ may be computed in $O(|V|2^{|V|-1})$ additions via
the Fast M\"{o}bius Transform \citep{kennes:mobius:91}.  ADtrees
\citep{moore:1998} provide a memory-efficient data-structure for storing
M\"{o}bius parameters.  Moreover, the matrix for the map $\nu$ can be shown
to have a Kronecker product structure; compare \citet{jokinen:2006}.

\begin{example}\rm({\it Two binary random variables}). \label{ex:two}
  \rm 
  Consider two binary random variables, i.e., $V=\{1,2\}$. Then the
  M\"{o}bius parameters are
  \[
    q_1 = p_{00}+p_{01},\quad
    q_2 = p_{00}+p_{10},\quad
    q_{12} = p_{00}.
  \]
  The joint cell probabilities can be recovered as
  \begin{align*}
    p_{00} &= q_{12}, &  p_{01} &= q_{1}-q_{12},\\
    p_{10} &= q_2 -q_{12}, &  p_{11} &= 1-q_1-q_2+q_{12}.
  \end{align*}
  The M\"{o}bius simplex is defined by the linear equalities expressing 
  that $p_i$, written in terms of $q$, is in the unit interval $[0,1]$ for
  all $i\in\cI$. In this example
  \[
  Q=\big\{ q=(q_1,q_2,q_{12})\in [0,1]^3 \,:\, q_1+q_2-1 \le q_{12} \le
  \min\{q_1,q_2\}   \big\}
  \]
  is a 3-dimensional simplex with vertices $(0,0,0)^t$, $(0,1,0)^t$,
  $(1,0,0)^t$, $(1,1,1)^t$.
\end{example}

As we show next, the constraints defining the independence model $\bP(\G)$
take on a simple form when expressed in terms of the M\"obius parameter
coordinates.  
\begin{theorem}
  \label{thm:moebconstraints}
  A probability vector $p\in\Delta$ belongs to the binary bi-directed graph
  model $\bP(\G)$ if and only if its M\"obius parameters $q=\mu(p)$ satisfy
  that for every disconnected set $D\subseteq V$,
  \begin{equation}\label{eq:disconnect}
    q_{D} = q_{C_1}q_{C_2}\cdots q_{C_r},
  \end{equation}
  where $C_1,\ldots,C_r$ are the inclusion-maximal connected sets forming
  the partition (\ref{eq:partition}).
\end{theorem}
\begin{proof}
  By Lemma \ref{lem:conncheck}, $p\in\bP(\G)$ implies
  (\ref{eq:disconnect}).  Conversely, consider a vector $q\in Q$ satisfying
  (\ref{eq:disconnect}), and let $p=\nu(q)$ be the associated probability
  vector. We show that $p\in\bP(\G)$ by verifying condition
  (\ref{eq:product}) in Lemma \ref{lem:conncheck}.  We proceed by induction
  on the number of ones in the vector $i_{D}\in\{0,1\}^D$ appearing in
  (\ref{eq:product}), for some $D$, which we denote by $k\in\{0,1,\ldots,V\}$.
  
  By (\ref{eq:disconnect}), the claim (\ref{eq:disconnect}) holds for
  $k=0$. Suppose that the claim holds for all $j<k$. Let $v$ be such that
  $i_{v}=1$ in $i_{D}$.  Let $C_{1},\ldots C_{r}$ be the partition of
  $D$ into inclusion-maximal connected components, and suppose that $v \in
  C_{\ell}$. Then
  \begin{eqnarray*}
    P(X_{D}=i_{D})&=& P(X_{D\setminus \{v\}}=i_{D\setminus \{v\}})
    - P(X_{D\setminus \{v\}}=i_{D\setminus \{v\}}, X_{v}=0)\\
    &=& \left[ P(X_{C_{\ell}\setminus \{v\}}=i_{C_{\ell}\setminus \{v\}})
      - P(X_{C_{\ell}\setminus \{v\}}=i_{C_{\ell}\setminus \{v\}}, X_{v}=0)
    \right]\\ \nonumber
    &&\quad\quad\times \prod_{j\neq \ell} P(X_{C_{j}}=i_{C_{j}})\\
    &=& \prod_{j=1}^{r} P(X_{C_{j}}=i_{C_{j}}).
  \end{eqnarray*}
  The second equality follows from the induction hypothesis applied to
  $i_{D\setminus \{v\}}$, and to $\bar \imath_{D}=(i_{D\setminus \{v\}},0)$
  since both vectors contain less than $k$ ones. Hence, we have shown that
  (\ref{eq:product}) holds true for all disconnected sets $D\subseteq V$.
\end{proof}

\begin{example}\rm({\it Four cycle}).
  For the bi-directed graph in Figure \ref{fig:bg}(a) we have 13 M\"obius
  parameters associated with connected sets
  \[
  q_1,q_2,q_3,q_4,\quad q_{12}, q_{13}, q_{24},q_{34},\quad q_{123},
  q_{124}, q_{134}, q_{234},\quad q_{1234}.
  \]
  In order to define a distribution obeying $X_1\ind X_4$ and $X_2\ind
  X_3$, the M\"obius parameters of the two disconnected sets must satisfy
   $q_{14}=q_1q_4$ and $q_{23}=q_2q_3$.
\end{example}

Theorem \ref{thm:moebconstraints} can be read as providing a model
parametrization.  Let $Q_\G=\mu(\bP(G))$ be the M\"obius parameter vectors
defining a distribution in $\bP(G)$.  Let $\cC(\G )$ be the family of
non-empty connected sets of $\G$.  Define $T_\G$ to be the set of vectors
$(q_C\mid C\in\cC(\G ))\in\RRR^{\cC(\G )}$ of M\"obius parameters of
connected sets for which there exists a vector $\bar q\in Q_\G$ with $\bar
q_C=q_C$ for all $C\in \cC (\G )$.
\begin{corollary}
  \label{cor:connpara}
  Let $\nu_\G:T_\G\to \bP(\G)$ be the multilinear map defined by setting
  M\"obius parameters of disconnected sets equal to the expression in
  (\ref{eq:disconnect}), obtaining a vector $q\in\RRR^{2^V-1}$, and setting
  $p=\nu(q)\in\bP(\G)$.  Then $\nu_G$ is a bijection, and we call it the
  M\"obius parametrization of the model $\bP(\G)$.
\end{corollary}
Since the M\"obius parameters are related via inequalities (compare Example
\ref{ex:two}), this parametrization is not variation independent, but
nevertheless is useful.  The definition of the M\"obius parameters is clearly
not symmetric under re-labelling of the two states taken by the random
variables.  However, such re-labelling does not change the model $\bP(\G)$
because it is defined purely in terms of independence relations.

\begin{figure}[t]
  \centering
  \vspace{-1cm}
  \includegraphics[width=9.5cm]{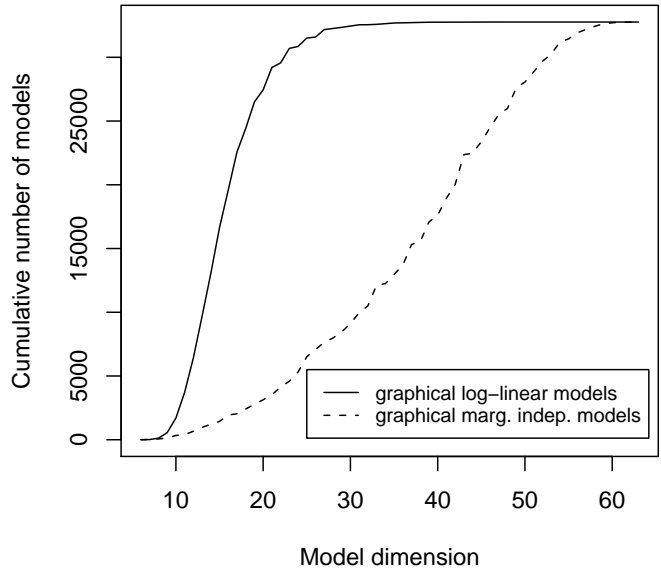} 
  \vspace{-0.5cm}
  \caption{Cumulative number of models per dimension for 
    $|V|=6$ binary variables.}
  \label{fig:dimplot}
\end{figure}

\begin{corollary}
  \label{cor:dim}
  The dimension of the model $\bP(\G)$ equals $\dim(\bP(\G))=|\cC(\G )|$,
  the number of non-empty connected sets in $\G$.
\end{corollary}
In contrast, the dimension of the (binary) graphical log-linear model
based on the undirected graph with the same edges as $\G$ would be
equal to the number of non-empty complete sets in $\G$.  Here a set
$A\subseteq V$ is complete if any two vertices in $A$ are adjacent.
Since every complete set is connected, the dimension of the model
$\bP(\G)$ is always larger than or equal to the dimension of the
corresponding graphical log-linear model; compare Figure
\ref{fig:dimplot}.

\begin{corollary}
  The family
  \[
  \bP_+(\G)= \{p\in \bP(\G) \,:\, p_i>0\; \mbox{ for all } i \in \cI\}  
  \]
  of distributions with positive joint cell probabilities in the binary
  bi-directed graph model $\bP(\G)$ forms a $|\cC(\G)|$-dimensional curved
  exponential family.
\end{corollary}
\begin{proof}
  More precisely stated, we claim that $\bP_+(\G)$ is a
  $|\cC(\G)|$-dimensional smooth manifold in the natural parameter space of
  the exponential family formed by the interior of the probability simplex
  $\Delta^o$.  Let $Q^o$ be the interior of the
  M\"obius simplex $Q$, and $Q_G^o$ the set of vectors in $Q^o$
  that satisfy the constraints (\ref{eq:disconnect}) in Theorem
  \ref{thm:moebconstraints}.  Let $d=2^{V}-|\cC(\G)|-1$ be the number of
  non-empty disconnected sets of $\G$.  Define the map
  $h:Q^o\to\RRR^{d}$ with coordinate functions
  $h_D(q)=q_D-q_{C_1}q_{C_2}\dots q_{C_r}$, where $C_1,C_2,\ldots,C_r$ form
  the inclusion-maximal connected set partition of the non-empty
  disconnected set $D$; compare (\ref{eq:partition}).  Since $h$ is
  $C^\infty$, it is clear that $Q_G^o=h^{-1}(0)$ is a $|\cC(\G)|$-dimensional
  smooth manifold in $\RRR^{2^V-1}$; compare e.g.~Thm.~1 in
  \cite{geiger:01}.  Our claim is now established because the
  diffeomorphism $\nu$ maps $Q^o$ to $\Delta^o$, and it
  is well-known that there is a diffeomorphism between $\Delta^o$
  (mean parameters) and the log-linear parameters (natural parameters of
  the exponential family).
\end{proof}

\begin{remark}\label{rem:dep-ratio}
  \rm Instead of using the M\"obius parameters in
  Theorem~\ref{thm:moebconstraints}, we could have employed the {\it dependence
  ratios}
  \[
  \tau_A = \frac{q_A}{\prod_{i\in A} q_i}
  \]
  introduced by \cite{ekholm:95}; see also \cite{ekholm:00,ekholm:03}, and
  \cite{darroch:83} where such ratios occur in specifying models termed
  {\it Lancaster additive}.  The ratio $\tau_A$ compares the probability
  $q_A$ computed from the joint distribution $p$ to the corresponding
  probability under the complete independence distribution that has the
  same univariate marginals as $p$.  Clearly,
  Theorem~\ref{thm:moebconstraints} also holds if we replace each M\"obius
  parameter by the corresponding dependence ratio.
\end{remark}

\section{Maximum likelihood estimation}
\label{sec:mle}

Assume we observe a sample of size $n$ drawn from a distribution $p$ in the
binary bi-directed graph model $\bP(\G)$, giving rise to multinomially
distributed counts $N(i)$, $i\in \cI$. \citep[For the link to Poisson
sampling see][ \S4.2.1.]{lau:bk} The probability of observing the
particular counts $n(i)\in\NNN_0$, $i\in \cI$, is equal to
\begin{equation}
  \label{eq:multinom}
  P(N(i)=n(i), \; i\in \cI)= \frac{n!}{\prod_{i\in\cI} n(i)!}
  \prod_{i\in\cI} p_i^{n(i)},
\end{equation}
where we set $0^0:=1$. Hence, the likelihood function for the model
$\bP(\G)$ is the map
\begin{equation}
  \label{eq:lik}
  \begin{split}
    L : \bP(\G) &\to \RRR,\\
    p &\mapsto \frac{n!}{\prod_{i\in\cI} n(i)!}
  \prod_{i\in\cI} p_i^{n(i)}.
  \end{split}
\end{equation}

\begin{proposition}
  \label{prop:mleexist}
  An MLE of $p\in\bP(\G)$ always exists.
\end{proposition}
\begin{proof}
  As a subset of the probability simplex $\Delta$, the model $\bP(\G)$
  is bounded. It is also closed, hence compact, which in conjunction
  with the continuity of the likelihood function implies the claim.
  Closedness follows from the fact that if for two sets $A,B\subseteq
  V$, $X_A\ind X_B$ under a sequence of probability distribution $P_n$
  with vector of joint cell probabilities $p_n\in \Delta$, then under
  a probability distribution $P$ corresponding to a limit point $p\in
  \Delta$ of the sequence $(p_n)$ it is also true that $X_A\ind X_B$;
  compare \citet[Prop.~3.12]{lau:bk}.
\end{proof}

If all counts $n(i)$, $i\in \cI$, are positive, then an MLE of
$p\in\bP(\G)$ will actually have positive joint cell probabilities, i.e.,
lie in $\bP_+(\G)$.  An open question is when an MLE exists in $\bP_+(\G)$
if some of the counts $n(i)$ are zero. For recent work on the analogous
question in the case of hierarchical log-linear models see
\cite{eriksson:2004}.  Another open problem concerns uniqueness of the MLE,
i.e., can one find a graph $\G$ and (non-degenerate) counts $n(i)$ such
that the likelihood function of $\bP(G)$ has more than one local maximum?

Ignoring an additive constant the log-likelihood function for the model
$\bP(\G)$ is of the form $\ell(p)= \sum_{i\in\cI} n(i)\log p_i$.
Using Proposition \ref{prop:satpara}, we can express the log-likelihood
function also in terms of M\"obius parameters as
\begin{equation*}
  \begin{split}
    \ell: Q_\G &\to \RRR,\\
    q&\mapsto \sum_{A\subseteq V} n(0_A,1_{V\setminus A})\log \bigg[
    \sum_{B: A \subseteq B} (-1)^{|B\setminus A|} q_B\bigg],
  \end{split}
\end{equation*}
where $q_\emptyset=1$.  Further, $\ell(q)$ can be written in terms of the
connected set M\"obius parameters $(q_C\mid C\in\cC(\G ))$ by replacing
$q_B$ for a disconnected set $B$ by the appropriate product of connected
set M\"obius parameters; see (\ref{eq:disconnect}).  

For two subsets $A,W\subseteq V$, nested as $A\subseteq W$, define
\[
p^W_A = P(X_{A}=0,X_{W\setminus A}=1).
\]
In particular, if $W=V$, then $p^{V}_A = P(X_{A}=0, X_{V\setminus A}=1)$ is
a joint cell probability. Similarly define $n_A^V$ to be the frequency of
observations in which $X_A = 0$, and $X_{V\setminus A} = 1$.  Then the
likelihood equations associated with the model $\bP(\G)$ are
\[
{{\partial \ell}\over {\partial q_C}} =
\sum_{A: \spo (C) \cap (A\setminus C) = \emptyset}
(-1)^{|C\setminus A|} {n_A^V \over p_A^V} p_{A\setminus C}^{V\setminus \spo(C)} = 0
\]
for every (non-empty) connected set $C$ in $\G$.  We prove this in the
Appendix (Corollary \ref{cor:likelihoodeqns}), where we also compute the
second derivative of $\ell(q)$, which yields the Fisher-information for
$\bP(G)$.

Having written the log-likelihood function as a function of the parameters
$(q_C\mid C\in\cC(\G))$, it can be maximized using gradient-based ascent
methods \citep[see also][]{lang:94,bergsma:rapcsak:2005}.  We implemented
such a method in the statistical programming environment R \citep{R} using
the routine `nlm'.  In doing this we found it beneficial to work with the
logarithms of the parameters $q_C$ because this linearizes
(\ref{eq:disconnect}); the examples we considered involve positive counts
such that we may assume that $q_C$ is positive and $\log q_C$ well-defined.
In our experience, this approach works well for smaller and sparser graphs
that induce a lower-dimensional model.  However, for larger and denser
graphs, such as in Figure~\ref{fig:trustgraph}(a), we found an alternative
approach that focuses on the model-defining constraints to perform better.
This alternative method, described in the next section, is
the binary analogue to the Iterative Conditional Fitting (ICF) algorithm
that was developed for ML fitting of Gaussian marginal independence models
\citep{drton:2003b,drton:2004jrssb}.  Binary ICF plays a role dual to the
Iterative Proportional Fitting (IPF) algorithm used to fit hierarchical
log-linear models.  

\section{Iterative conditional fitting}
\label{sec:icf}

Starting from some feasible estimate in $\bP(G)$, such as the uniform
distribution, the ICF algorithm improves a current feasible estimate by
cycling through the vertex set $V$ and performing an update step for each
one of the vertices. At the update step for variable $v\in V$ the marginal
distribution $P^{X_{-v}}$ of the variables $-v=V\setminus\{v\}$ is fixed,
and the conditional distribution $P^{X_v\mid X_{-v}}$ required to determine
the joint distribution of $(X_v\mid v\in V)$ is estimated.  This estimation
is done subject to constraints that ensure that the newly determined joint
distribution remains in the model $\bP(G)$.  In this presentation of ICF we
assume that all observed counts $n(i)$, $i\in\cI$, are positive, which in
particular entails that they were drawn from a distribution
$p\in\bP_+(\G)$. Moreover, maximizing the likelihood function over
$\bP(\G)$ is equivalent to maximizing it over the submodel $\bP_+(\G)$, and
we can assume that all joint distributions $P$ considered in the sequel
have positive joint cell probabilities $p_i>0$.  In the case of zero
counts, which will be considered in future work, the conditional likelihood
function considered in Algorithm 5.1 is still concave but need no longer be
strictly concave.  Hence, the possibility of optima on the boundary has to
be taken into account.
 
For fixed marginal probability $P(X_{-v}= i_{-v})$, the joint cell
probability $P(X_v = i_v, X_{-v}= i_{-v})=0$, $i\in \cI$, is determined by
the conditional parameter
\[
\theta_v(i_{-v}) = \theta_v(X_{-v}=i_{-v}) := P(X_v = 0 \mid X_{-v} =
i_{-v}).
\]
Let $\cI_{-v}=\{0,1\}^{V-1}$. Then there are $|\cI_{-v}|= 2^{V-1}$ many
parameters $\theta_v(i_{-v})$.  Notice that if $v\in D$ then
\begin{align}
  \nonumber
  P(X_D=0) &= \sum_{B: D\subseteq B} P(X_B=0, X_{V\setminus B}=1)\\
  \label{eq:construction}
  &= \sum_{
i_{-v} = (0_{B\setminus \{v\}},1_{V\setminus B}) \in \cI_{-v} \,:\,
D\subseteq B
} \theta_v ( i_{-v})
P(X_{-v}=i_{-v} ).
\end{align}

In general, the binary bi-directed graph model $\bP(\G)$ imposes
constraints on the conditional distribution $P^{X_v\mid X_{-v}}$.  In order
to specify the constraints in a non-redundant way, we focus on constraints
of the form (\ref{eq:disconnect}), rather than the equivalent conditional
independence restrictions.  Specifically, suppose that $D$ is a
disconnected set and that $C$ is the inclusion-maximal connected subset of
$D$ containing $v$. By equation (\ref{eq:disconnect}) we require
\begin{equation}\label{disconnect2}
P(X_D = 0) = P(X_C=0) P(X_{D\setminus C} = 0).
\end{equation}
Note that $D\setminus C$ may not be connected, so the model may require
further factorization of $P(X_{D\setminus C} = 0)$. However, this only
imposes a constraint on the fixed $P(X_{-v})$ margin, and so does not
concern us here.  We now express the constraint (\ref{disconnect2}) as
\begin{eqnarray}
\nonumber \lefteqn{P(X_v = 0 \mid X_{D\setminus\{v\}} = 0) P(X_{D\setminus\{v\}}
= 0)}\\
\label{disconnect3}& & = P(X_v = 0 \mid X_{C\setminus\{v\}}=0)
P(X_{C\setminus\{v\}} = 0) P(X_{D\setminus C} = 0).
\end{eqnarray}
(It is implicit here that if $C\setminus\{v\} = \emptyset$ then the second
term on the right hand side is omitted.) Observe that only the first terms
on each side depend on $\theta_v(\cdot )$. Using (\ref{eq:construction}),
the first term on the left hand side of (\ref{disconnect3}) may be
expressed as
\begin{eqnarray}
\nonumber
\displaystyle
\lefteqn{P(X_v = 0 \mid X_{D\setminus\{v\}}=0) }\\
\nonumber
 &=& \sum_{
j \in \{0,1\}^{V\setminus D} }
 P(X_v=0, X_{V\setminus D} = j
\mid  X_{D\setminus \{v\}} = 0).\\
\label{eq:lin1}
 &=& \sum_{
j \in \{0,1\}^{V\setminus D} }
\theta_v(X_{V\setminus D} = j, X_{D\setminus \{v\}} = 0)
 P(X_{V\setminus D} = j \mid  X_{D\setminus \{v\}} = 0).
\end{eqnarray}
Similarly, the first term on the right hand side of (\ref{disconnect3}) may
be expressed as
\begin{eqnarray} 
\nonumber
\displaystyle
\lefteqn{P(X_v = 0 \mid X_{C\setminus\{v\}}=0) }\\
\label{eq:lin2}
 &=& \sum_{j \in \{0,1\}^{V\setminus C}}
\theta_v(X_{V\setminus C} = j, X_{C\setminus \{v\}}
= 0 ) P(X_{V\setminus C} = j
\mid  X_{C\setminus \{v\}} = 0).
\end{eqnarray}

Now, if the set $D\setminus \{v\}$ was connected, then $C$ and $D\setminus C$
would also be connected, contrary to the assumption. Since in the ICF
algorithm we assume that all constraints on the marginal distribution
of $X_{V\setminus \{v\}}$ hold, it follows that
\[
P(X_{D\setminus \{v\}} = 0) = P(X_{D\setminus C} = 0) P(X_{C\setminus \{v\}} = 0).
\]
(Again, both $C\setminus\{v\}$ and $D\setminus C$ may not be connected, so these
terms may factorize further.) Since these terms are non-zero, they
cancel from both sides of (\ref{disconnect3}), leaving the constraint
\begin{multline}
 \displaystyle\sum_{
j \in \{0,1\}^{V\setminus D} }
\theta_v(j, 0_{D\setminus\{v\}} )  P(X_{V\setminus D} =
j \mid  X_{D\setminus \{v\}} = 0)
\label{eq:lin3}
 = \\
 \sum_{
j \in \{0,1\}^{V\setminus C} 
}
\theta_v(j, 0_{C\setminus\{v\}} )  P(X_{V\setminus C} =
j
\mid  X_{C\setminus \{v\}} = 0).
\end{multline}
It is important to note that for fixed margin $P^{X_{-v}}$ the constraints
(\ref{eq:lin3}) are linear in the conditional parameters $\theta_v$.
The full set of constraints on the $\theta_v$ parameters may be obtained by
considering every disconnected set $D$ containing $v$ and identifying the
inclusion-maximal connected set $C\subset D$ containing $v$. 

Let
\[
\mathfrak{D}_v = \left\{ D \,:\, D\subseteq V,\;  v\in D,\; D \hbox{ is
disconnected} \right\}.
\]
For each set $D\in \mathfrak{D}_v$, we define $C_v(D)$ to be the
inclusion-maximal connected subset of $D$ containing $v$.  The disconnected
sets $\mathfrak{D}_v$ and the connected sets $C_v(D)$ can be computed in
preprocessing.  Then the ICF update for vertex $v$ can be implemented as
follows.

\begin{algorithm}{\it Update step in Iterative Conditional Fitting.}
\label{alg:icf}
\rm \mbox{}\\
{\sl Input:} A probability vector $p\in\bP(G)$ and vertex $v$.\\
{\sl Output:} A probability vector $\bar p\in\bP(G)$ such that $L(\bar
p)\ge L(p)$.\\
{\sl Step 1.} Construct the $\mathfrak{D}_v\times \cI_{-v}$ constraint
matrix $A = (a_{rs})$, where for each pair $(D_r,j_s)\in \mathfrak{D}_v\times
\cI_{-v}$ we set 
\begin{eqnarray*}
  \lefteqn{a_{rs} =  P(X_{V\setminus D_r}=(j_s)_{V\setminus D_r}
    \mid X_{D_r\setminus \{v\}}=0) \,
    \hbox{I}\{(j_s)_{D_r\setminus \{v\}}=0\}}\\[6pt]
  & & \quad\quad\quad - P(X_{V\setminus C_v(D_r)}=(j_s)_{V\setminus C_v(D_r)}
  \mid X_{C_v(D_r)\setminus \{v\}}=0) \,
  \hbox{I}\{(j_s)_{C_v(D_r)\setminus \{v\}}=0\}. 
\end{eqnarray*}
Here all probabilities are computed under the distribution induced by the
probability vector $p$, and $\hbox{I}(\cdot )$ is the indicator
function.\\
{\sl Step 2.} Maximize the strictly concave conditional log-likelihood
function
\[
 \sum_{i_{-v}\in\cI_{-v}} n(i_{-v},0) \log \theta(i_{-v}) +
n(i_{-v},1) \log \{1-\theta(i_{-v})\}
\]
subject to the linear constraints $A\theta = 0$, where
$\theta=(\theta(i_{-v})\mid i_{-v}\in\cI_{-v})$ is the vector of all
conditional parameters. (If all counts are positive, the inequality constraints
$\theta\in [0,1]^{\cI_{-v}}$ need not be considered explicitly.)\\
{\sl Step 3.}  Use the solution $\theta_v$ from step 2 to compute the new
probability vector $\bar p\in\bP(G)$ via
\[
\bar p_{i_vi_{-v}}=
\bar P(X_v = i_v, X_{-v} = i_{-v}) = 
\begin{cases}
  \theta_v(i_v)\, P(X_{-v} = i_v) &\hbox{if } 
  i_v=0,\\ 
  {}[1-\theta_v(i_v) ]\, P(X_{-v} = i_{-v}) &\hbox{if } 
  i_v=1.
\end{cases}
\]
\end{algorithm}

The optimization problem in step 2 of the ICF update algorithm has a unique
local maximum and is not difficult to solve. For example, one can employ
the gradient projection method \cite[ \S2.3]{bertsekas:1999}, which
performs a line search along the direction of the gradient projected on the
kernel of $A$.  A line search based on the Armijo-rule ensures convergence
of the gradient projection method.  The computation of the $2^{V-1}\times
2^{V-1}$ projection matrix $I-A'(AA')^{-1}A$ requires the inversion of the
$\mathfrak{D}_v\times \mathfrak{D}_v$ matrix $AA'$ which is of full rank.
However, the projection matrix has to be computed only once in order to
solve the optimization problem in step (4) of Algorithm \ref{alg:icf}.
Since the Hessian of the conditional log-likelihood function maximized in
step 2 of Algorithm \ref{alg:icf} is diagonal it is also feasible to employ
second derivative information in a projected Newton method, in which
$\theta$ is scaled by the matrix with diagonal elements equal to one over
the square root of the diagonal elements of the Hessian.  Since the Hessian
depends on $\theta$, the projection matrix in a projected Newton method has
to be recomputed every time $\theta$ is updated. However, based on our
experience with our implementation of ICF in R, employing the Hessian
information is beneficial.

Having tackled the individual ICF updates we can run ICF from a feasible
starting value.  The algorithm then produces a sequence of feasible
estimates whose accumulation points are solutions to the likelihood
equations.  In fact, the sequence is guaranteed to converge if there exist
only finitely many solutions to the likelihood equations.  These
convergence guarantees follow from general results about iterative partial
maximization algorithms \citep[Appendix]{drtoneichler:2006}, of which ICF
is an incarnation.  Here, `partial maximization' refers to the fact that
the update step for vertex $v$ in Algorithm \ref{alg:icf} maximizes the
log-likelihood function $\ell(q)$, $q=(q_C\mid C\in\cC(G))$ partially,
namely when varying only components $q_C$ for which $C$ is a connected set
containing vertex $v$.  Components $q_C$ for connected sets not containing
$v$ remain fixed at the current estimates.

In the above we proceeded vertex-by-vertex and estimated the univariate
conditional distribution of $X_v$ given $X_{-v}$.  In the Gaussian case,
\cite{drton:2004jrssb} describe how to run the ICF algorithm with
multivariate updates.  In this variant, one chooses complete vertex sets
$C\subseteq V$ and estimates, for fixed margin of $X_{-C}=X_{V\setminus
  C}$, the (multivariate) conditional distribution of $X_C$ given $X_{-C}$
under the marginal independence constraints.  Such multivariate updates are
also possible in the binary case discussed here.  Let 
\[
\theta_C(i_{-C}) = \theta_C(X_{-C}=i_{-C}) := P(X_C = 0 \mid X_{-C} =
i_{-C}).
\]
Let $\cI_{-C}=\{0,1\}^{V-|C|}$. 
Then (\ref{eq:construction}) becomes
\begin{align}
  \label{eq:constructionmulti}
  P(X_D=0) &= 
   \sum_{ i_{-C} = (0_{B\setminus \{C\}},1_{V\setminus B}) \in \cI_{-C}
    \,:\, D\subseteq B } \theta_C ( i_{-C}) P(X_{-C}=i_{-C} ).
\end{align}
Since the set $C$ is complete there are no equality constraints among the M\"obius
parameters $q_A$, $\emptyset\not=A\subseteq C$, and one can proceed
similarly as in the discussion leading up to (\ref{eq:lin3}) to devise an
analog to Algorithm \ref{alg:icf} with multivariate updates over complete
sets.

\section{Example: Social survey data}
\label{sec:survey}

Sociologists and political theorists have long been interested in the
relationship between trust in social institutions and trust in other
members of society \citep{putnam:bowling:2002, sztompka:trust:2000,
  levi:trust:1998}. Here, as an illustration of an exploratory analysis
using binary independence models we examine seven questions relating to
trust that are taken from the U.S.~General Social Survey during the years
1975-94:

\begin{itemize}
\item{\sc Trust}\par
   {\it Generally speaking, would you say that most people can 
be trusted or that you can't be too careful in life.} (Can Trust; Cannot Trust; Depends)
 \item{\sc Helpful}\par
{\it Would you say that most of the time
people try to be helpful, or that they are mostly just looking out
for themselves?} (Helpful; Lookout for Self; Depends)
\item{\sc Confidence in institutions}\par
{\it I am going to name some institutions in this country. As 
far as the people running these institutions are concerned, 
would you say you have a great deal of confidence, only some 
confidence, or hardly any confidence at all in them?} (A Great Deal; Only Some; Hardly Any)
\begin{itemize}
\item[]{\sc ConClerg:} Organized religion
\item[]{\sc ConLegis:} Congress
\item[]{\sc ConBus:} Major Companies
 \end{itemize}
\item {\sc Membership of organizations}\par
{\it Here is a list of various organizations. Could you tell me whether or not you are 
a member of each type?} (Yes; No)
\begin{itemize}
\item[]{\sc MemUnion:} Labour unions
\quad\quad\quad{\sc MemChurch:} Churches
\end{itemize}
\end{itemize}
There were $13,486$ individuals who gave valid responses to all of these
questions.  For the purposes of illustration, for the questions relating to
confidence in institutions we combine `Some' and `Hardly Any' to form a
`No' response; similarly for {\sc Trust} we combine `Cannot trust' with
`Depends' to form a `No' group, and for {\sc Helpful} we combine `Take
advantage' with `Depends' to form a `No' group. The counts are displayed in
Table \ref{tab:datatrust}.

\begin{table}
 \caption{Data from the U.S.~General Social Survey relating to
   Trust.\label{tab:datatrust}}  
  \centering 
  \begin{tabular}{lllll  rrrr}
    \hline
 &&&&\multicolumn{1}{c}{ }&\multicolumn{4}{c} {\sc\small Helpful}\\
                        &&&&                      & \multicolumn{1}{c}{Yes}&\multicolumn{1}{c}{Yes} & \multicolumn{1}{c}{No}&
                        \multicolumn{1}{c}{No}\\      
                         \cline{6-9}                           
 {\sc\small Con.} & {\sc\small Con.} & {\sc\small Con.} & {\sc\small Mem.} & {\sc\small Mem.}
 & \multicolumn{2}{c} {\sc\small Trust} & \multicolumn{2}{c} {\sc\small Trust}\\   
 {\sc\small Bus.} & {\sc\small Clerg.} & {\sc\small Legis.} & {\sc\small Church} & {\sc\small Union} 
 &      \multicolumn{1}{c}{Yes} &   \multicolumn{1}{c}{No} &   \multicolumn{1}{c}{Yes} &   \multicolumn{1}{c}{No}\\  
\hline\hline                          
Yes   &   Yes  &      Yes  &     Yes  &      Yes  &    18  &  4  &  5  &  5\\
 Yes      &     Yes  &       Yes &    Yes     &      No   &    79 &  47 &  17 &  30\\
 \cline{4-9}
  Yes&           Yes   &         Yes  &   No  &    Yes    &  8   & 9   & 1  & 15\\
    Yes&           Yes &        Yes  &   No     &     No    & 88 & 55 &  22 &  79\\
     \cline{3-9}
   Yes&             Yes&        No  &     Yes    &    Yes   &  22   & 11   & 10  & 13\\
   Yes&              Yes&            No   &        Yes       &   No   & 194   & 95   & 33   & 77\\
    \cline{4-9}
    Yes&              Yes&         No      &    No       &  Yes  &  31  & 10   & 13   & 23\\
    Yes&               Yes&        No   &     No     &  No   & 179   & 82   & 58  & 122\\
    \cline{2-9}
    Yes&              No  &   Yes    &       Yes    &  Yes   &     7    & 5    & 1   & 3\\
     Yes&               No&           Yes &         Yes    &  No  & 40          & 27    & 11   & 23\\
      \cline{4-9}
     Yes&                No&          Yes &    No    &  Yes  &  9        & 10    & 1   & 12\\
      Yes&                No&          Yes   &           No    &  No   & 68         & 56    & 33   & 73\\
       \cline{3-9}
      Yes&                   No       &     No   &    Yes   &  Yes   &  15      & 13    & 6   & 14\\
      Yes&                      No     &   No  &      Yes   &   No   & 188        & 117   & 52  & 100\\
       \cline{4-9}
       Yes&                  No     &        No&    No    &    Yes &  32        & 28    & 22   & 35\\
       Yes&                    No&                  No&      No    &    No & 366        & 185   & 120  & 312\\
       \hline 
         No&                   Yes&       Yes&  Yes    &   Yes &     7     & 5    & 2    & 6\\
          No&                   Yes&                  Yes&   Yes     &  No   & 62           & 32       & 11   & 48\\
           \cline{4-9}
         No&                Yes&                   Yes&   No    & Yes  & 5           & 9        & 2   & 12\\
         No&                  Yes&                       Yes&   No    &  No  & 38           & 37       & 11   & 64\\
          \cline{3-9}
          No&           Yes   &                   No&   Yes   &  Yes &  40        & 26    & 17   & 34\\
          No&                  Yes&                     No&   Yes   & No  & 270         & 187    & 73  & 281\\
           \cline{4-9}
          No&                 Yes&                        No&   No  & Yes  &  25         & 33     & 11   & 50\\
           No&                  Yes&                         No&  No  & No   & 202          & 216    & 84  & 356\\
           \cline{2-9}
       No&    No &                     Yes&  Yes  & Yes   &    5        & 2     & 3   & 11\\
         No&                     No&                       Yes& Yes  & No   & 51           & 32  & 17   & 59\\
          \cline{4-9}
         No&                   No&                  Yes&  No  & Yes  &  15          & 18   & 7   & 33\\
         No&                      No&                    Yes&  No  & No    & 104          & 79    & 40  & 172\\
          \cline{3-9}
         No&       No&                No&   Yes & Yes   &  74        & 62  & 27  & 108\\
          No&                  No&                     No&  Yes  & No  & 603          & 469   & 177  & 654\\
           \cline{4-9}
          No&                  No&                   No& No  & Yes   &  199         & 181   & 84  & 305\\
          No&                   No&                No& No  & No  & 1002         & 920  & 460 & 1818\\
\hline
  \end{tabular}
\end{table}

\begin{figure}[tbp]
  \vspace{1cm}
  \begin{center}
    \small
  \hspace{-1.75cm}
    \newcommand{\myNode}[2]{\circlenode{#1}{\makebox[3.1ex]{#2}}}
    \begin{pspicture}(0,0)(2,3)
    \psset{linewidth=0.6pt,arrowscale=1.25 2,arrowinset=0.1}  
    \psset{fillcolor=lightgray, fillstyle=solid}          
    \rput(-1,-.1){(a)}
    \rput(0, 0){\myNode{1}{\vbox{\sc\tiny Mem.\\ Union}}}
    \rput(2, 0){\myNode{2}{\vbox{\sc\tiny Con.\\Bus.}}} 
    \rput(0,1.75){\myNode{3}{\vbox{\sc\tiny Mem.\\Church}}} 
    \rput(4,1.75){\myNode{4}{\vbox{\sc\tiny Con.\\Clerg.}}} 
    \rput(4, 0){\myNode{5}{\vbox{\sc\tiny Con.\\Legis.}}}
    \rput(2,2.6){\myNode{6}{{\sc\tiny Trust}}} 
    \rput(3.5,3.25){\myNode{7}{\vbox{\sc\tiny Help-\\ful}}} 
     \psset{fillcolor=white, fillstyle=solid}    
    \ncline{<->}{2}{3}
     \ncline{<->}{2}{7}
     \psset{fillcolor=white, fillstyle=none} 
     \nccurve[ncurv=1,angleB=170,angleA=55]{<->}{3}{7}
     \ncline{<->}{3}{4}
     \ncline{<->}{7}{4}
     \ncline{<->}{4}{6}
      \ncline{<->}{6}{2}
     \ncline{<->}{2}{5}
     \ncline{<->}{2}{4}
     \ncline{<->}{4}{5}
     \ncline{<->}{6}{7}
     \ncline{<->}{3}{6}
      \psset{linestyle=dashed}  
  \ncline{<->}{1}{2}
   \nccurve[ncurv=1,angleB=240,angleA=120]{<->}{1}{3}
    \end{pspicture}
\hskip120pt
 \begin{pspicture}(0,0)(2,3)
    \psset{linewidth=0.6pt,arrowscale=1.25 2,arrowinset=0.1}  
    \psset{fillcolor=lightgray, fillstyle=solid}          
    \rput(-1,-.1){(b)}
    \rput(0, 0){\myNode{1}{\vbox{\sc\tiny Mem.\\ Union}}}
    \rput(2, 0){\myNode{2}{\vbox{\sc\tiny Con.\\Bus.}}} 
    \rput(0,1.75){\myNode{3}{\vbox{\sc\tiny Mem.\\Church}}} 
    \rput(4,1.75){\myNode{4}{\vbox{\sc\tiny Con.\\Clerg.}}} 
    \rput(4, 0){\myNode{5}{\vbox{\sc\tiny Con.\\Legis.}}}
    \rput(2,2.6){\myNode{6}{{\sc\tiny Trust}}} 
    \rput(3.5,3.25){\myNode{7}{\vbox{\sc\tiny Help-\\ful}}} 
     \psset{fillcolor=white, fillstyle=none} 
    \ncline{-}{2}{3}
    \ncline{-}{2}{1}
     \ncline{-}{2}{5}
      \ncline{-}{2}{4}
     \ncline{-}{2}{7}
     \ncline{-}{2}{6}
     \ncline{-}{1}{3}
     \ncline{-}{4}{3}
     \ncline{-}{1}{7}
      \ncline{-}{1}{4}
     \ncline{-}{3}{6}
     \nccurve[ncurv=1,angleB=170,angleA=55]{-}{3}{7}
     \ncline{-}{4}{5}
     \ncline{-}{4}{6}
     \ncline{-}{6}{7}
    \end{pspicture}
  \end{center}
\vskip5pt
  \caption{ \label{fig:trustgraph}Analysis of Trust data. 
  (a) Marginal independence model; dashed edges correspond to pairwise odds ratios less than one. (b) Classical graphical log-linear model.  }
\end{figure}
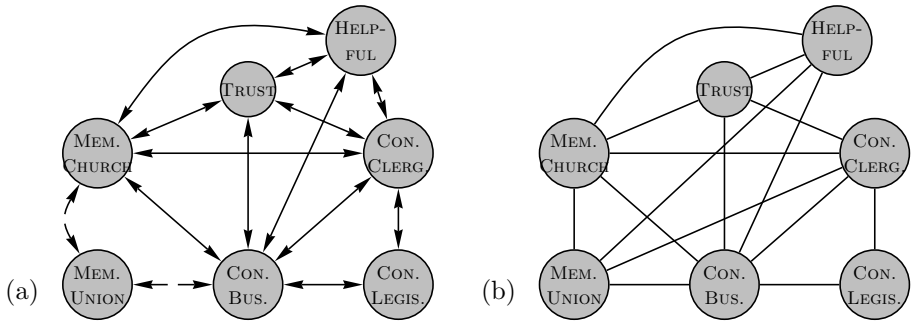
 
Using ICF in a backward stepwise selection we found the graph shown in
Figure \ref{fig:trustgraph}(a).  Assuming that the data in Table
\ref{tab:datatrust} arose in multinomial sampling, we obtain a deviance of
$32.67$ over $26$ degrees of freedom, when compared to the saturated model
of no independence.  Using an asymptotic $\chi^2$-approximation a p-value
of $0.172$ is obtained and the model is found not to be contradicted by the
data.  Since some expected cell counts are small, the asymptotic
approximation should be treated with some caution.  In the selected model all
variables are marginally associated with confidence in business, but it is
interesting that confidence in congress is marginally associated only to
the two other confidence variables.  Similarly, union membership is not
marginally associated with additional variables other than church membership
 and confidence in business; the graph implies 
\begin{align*}
    \hbox{\sc ConLegis}\; &\ind\; \hbox{\sc Helpful, Trust, MemUnion,
    MemChurch}, \quad\text{and}\\
  \hbox{\sc MemUnion}\; &\ind\; \hbox{\sc Helpful, Trust, ConLegis,
    ConClerg}.
\end{align*}
It is perhaps of little surprise that in the fitted distribution the
marginal odds ratio between {\sc MemUnion} and {\sc ConBus} is less than
one; it equals 0.83.  Except for the odds ratio between {\sc MemUnion} and
{\sc MemChurch}, which is equal to 0.85, all other fitted pairwise odds
ratios are greater than or equal to $1$.
 
For purposes of comparison, in Figure \ref{fig:trustgraph}(b) we include a
classical graphical log-linear model obtained using the MIM program
\citep{edwards:2000} by backward stepwise selection, among all undirected
models.  This model has a deviance of $87.62$ over $88$ degrees of freedom.
When comparing the undirected and bi-directed models it is quite striking
that the undirected model contains one more edge, yet $62$ fewer
parameters. Observe that in the undirected graph union membership is
also adjacent to the variables relating to confidence in clergy and whether
or not people are helpful.  We remark that latent variable models
could be used for further analyses of these data.

\section{Independence and symmetry}
\label{sec:symm-indep}

In this section we demonstrate how symmetry can be incorporated in the
marginal independence models proposed earlier.  The issue of symmetry
naturally arises for the twin data shown in Table \ref{tab:mzdata} in the
introduction.
  Recall that we observe four binary indicators which inform us about each
  twins' alcohol dependence ($A_i$) and depression status ($D_i$).  
  When inspecting Table \ref{tab:mzdata}, one
notices that counts related by exchanging the index labels $1$ and $2$ are
often very similar.

 \begin{table}
 \caption{\label{tab:symmle} ML estimate of the joint distribution 
       under the symmetry group  ${\mathcal S}_{\hbox{\scriptsize twin}}$;
empirical distribution in parenthesis.}
  \centering
  \begin{tabular}{ccccccc}
    \hline
    & & \multicolumn{2}{c}{$D_1=0$} & & \multicolumn{2}{c}{$D_1=1$}\\
    \cline{3-4} \cline{6-7}
    & & \multicolumn{1}{c}{$D_2=0$} & \multicolumn{1}{c}{$D_2=1$} & &
    \multicolumn{1}{c}{$D_2=0$} & \multicolumn{1}{c}{$D_2=1$}\\ 
    \hline
    \multirow{4}{*}{$A_1=0$} 
    &  \multirow{2}{*}{$A_2=0$}
                      & 0.4824 & 0.1441 && 0.1441 & 0.0854 \\
            &         &(0.4824)&(0.1340)&&(0.1541)&(0.0854)\\
            \cline{2-7}
    &  \multirow{2}{*}{$A_2=1$}
                      & 0.0193 & 0.0142 && 0.0092 & 0.0159\\
            &         &(0.0251)&(0.0151)&&(0.0117)&(0.0168)\\
            \hline
    \multirow{4}{*}{$A_1=1$}
    &  \multirow{2}{*}{$A_2=0$}
                      & 0.0193 & 0.0092 && 0.0142 & 0.0159\\
            &         &(0.0134)&(0.0067)&&(0.0134)&(0.0151)\\
            \cline{2-7}
    &  \multirow{2}{*}{$A_2=1$}
                      & 0.0050 & 0.0050 && 0.0050 & 0.0117\\
            &         &(0.0050)&(0.0034)&&(0.0067)&(0.0117)\\
            \hline
  \end{tabular}
\end{table}

Let $\mathcal{S}$ be a group of permutations on the index set $V$.
  The group $\mathcal{S}$ acts on the set of
elementary joint events $\cI$ by permuting the components of $i\in\cI$.  In
other words, for $\sigma\in\mathcal{S}$ and $i=(i_v\mid v\in V)\in\cI$, we
define $\sigma(i) = (i_{\sigma(v)}\mid v\in V)$.  This action induces the
symmetry model
\[
\bP(\mathcal{S})=\{p\in \Delta\mid
p_i=p_{\sigma(i)}\;\forall\sigma\in\mathcal{S}\}.
\]
\begin{example}\rm({\it Twin data}).
   The symmetry group  
  \begin{equation}\label{eq:stwin}
  \mathcal{S}_{\hbox{\scriptsize twin}}=\{(A_1)(A_2)(D_1)(D_2),\, (A_1\, A_2)(D_1\, D_2)\}
  \end{equation}
   represents symmetry when exchanging vertex $A_1$ with $A_2$, and at the same time exchanging $D_1$ with $D_2$.
  This symmetry corresponds to irrelevance of the labels given to the two twins.
  \end{example}
  
  Since the symmetry model is a linear exponential family, the MLE may be computed by 
  simply averaging the empirical cell counts over the orbit induced by the group action:
  \[
  \hat{p}_{\mathcal{S}}(i) = {1\over |S(i)|}  \sum_{j \in S(i)} {n(j)\over n}
  \]
  where $S(i) = \{ \sigma(i) \mid \sigma \in \mathcal{S}\}$ is the orbit of cell $i$ under
  the group $\mathcal{S}$,  $n(j)$ is the empircal count for cell $j$, and $n$ is the total sample size.
   For the twin data the ML estimate is shown in Table \ref{tab:symmle}. The deviance is $4.62$
   on $6$ degrees of freedom, indicating a good fit. We now turn to testing the marginal independence hypothesis mentioned in the introduction, in conjunction with symmetry.
 
 A permutation $\sigma\in\mathcal{S}$ induces a new graph $G_\sigma$ by
renaming vertex $v\in V$ to $\sigma(v)\in V$.  In other words, the graph
$G_\sigma$ has the same vertex set $V$ as the original graph $G=(V,E)$ but
there is an edge $v\bi w$ in $G_\sigma$ if and only if there is an edge
$\sigma^{-1}(v)\bi\sigma^{-1}(w)$ in the original graph $G$.
We say that a group of permutations $\mathcal{S}$ leaves the graph $G$ {\em
  invariant\/} if $G_\sigma=G$ for all $\sigma\in\mathcal{S}$, in other words, 
  $\mathcal{S}$ is a subgroup of the automorphism group of $G$.
   It follows
that no new independences are introduced when imposing symmetry on the
distributions in $\bP(G)$.  We will restrict attention to this case in what follows.
 
 \addtocounter{theorem}{-1}
  \begin{example}\rm({\it continued}).
  \rm Let $G$ be the graph displayed in Figure \ref{fig:bg}(a), under
  the variable-vertex correspondence $(1,2,3,4)=(A_1,A_2,D_1,D_2)$
  the independence pattern is $A_1\ind D_2$ and $A_2\ind D_1$. 
  The  group $\mathcal{S}_{\hbox{\scriptsize twin}}$ given in (\ref{eq:stwin}) leaves $G$ invariant.
  \end{example}
   \addtocounter{theorem}{4}

\begin{table}%
 \caption{ \label{tab:twinsymmle} ML estimate of the joint distribution under $A_1\ind D_2$ and $A_2\ind D_1$ together with
  the symmetry group $\mathcal{S}_{\hbox{\scriptsize twin}}$;
empirical distribution in parenthesis.}
  \centering
  \begin{tabular}{ccccccc}
    \hline
    & & \multicolumn{2}{c}{$D_1=0$} & & \multicolumn{2}{c}{$D_1=1$}\\
    \cline{3-4} \cline{6-7}
    & & \multicolumn{1}{c}{$D_2=0$} & \multicolumn{1}{c}{$D_2=1$} & &
    \multicolumn{1}{c}{$D_2=0$} & \multicolumn{1}{c}{$D_2=1$}\\ 
    \hline
    \multirow{4}{*}{$A_1=0$} 
    &  \multirow{2}{*}{$A_2=0$}
                      & 0.4612 & 0.1486 && 0.1486 & 0.0957 \\
            &         &(0.4824)&(0.1340)&&(0.1541)&(0.0854)\\
            \cline{2-7}
    &  \multirow{2}{*}{$A_2=1$}
                      & 0.0249 & 0.0204 && 0.0057 & 0.0104\\
            &         &(0.0251)&(0.0151)&&(0.0117)&(0.0168)\\
            \hline
    \multirow{4}{*}{$A_1=1$}
    &  \multirow{2}{*}{$A_2=0$}
                      & 0.0249 & 0.0057 && 0.0204 & 0.0104\\
            &         &(0.0134)&(0.0067)&&(0.0134)&(0.0151)\\
            \cline{2-7}
    &  \multirow{2}{*}{$A_2=1$}
                      & 0.0100 & 0.0038 && 0.0038 & 0.0054\\
            &         &(0.0050)&(0.0034)&&(0.0067)&(0.0117)\\
    \hline
  \end{tabular}
\end{table}


\begin{theorem}
  \label{thm:symm-indep-1}
  If the symmetry group $\mathcal{S}$ leaves the graph $G$ invariant, then
  a distribution $p\in\bP(G)$ is in the symmetry model $\bP(\mathcal{S})$
  if and only if the M\"obius parameters $q\in Q_G$ for $p$ satisfy that
  $q_C=q_{\sigma(C)}$ for all connected sets $C$ in $G$.
\end{theorem}
\begin{proof}
  First, note that under the assumed invariance of the graph, a set
  $C\subseteq V$ is connected in $G$ if and only if $\sigma(C)$ is
  connected for all $\sigma\in\mathcal{S}$.
  
  Consider $p\in\bP(G)\cap\bP(\mathcal{S})$, and let $C\in\cC(G)$ and
  $\sigma\in\mathcal{S}$.  Since $i_C=\sigma^{-1}(i)_{\sigma(C)}$ we obtain
  that
  \[
  q_C = \sum_{i\in\cI: i_C=0} p_i = 
  \sum_{i\in\cI: i_C=0} p_{\sigma^{-1}(i)} = 
  \sum_{j\in\cI: j_{\sigma(C)}=0} p_j = q_{\sigma(C)}.
  \]
  
  Conversely, assume that the M\"obius $q\in Q_G$ satisfy that
  $q_C=q_{\sigma(C)}$ for all $C\in\cC(G)$.  Let $D\subseteq V$ be
  disconnected and uniquely partitioned into inclusion-maximal connected
  sets as $D = C_1\dot\cup C_2\dot\cup \dots \dot\cup C_r$.
  Then the unique decomposition of $\sigma(D)$ into inclusion-maximal
  connected sets is given by
  \[
  \sigma(D)=\sigma(C_1)\dot\cup \sigma(C_2)\dot\cup \cdots \dot\cup
  \sigma(C_r),
  \]
  which implies that 
  \[
  q_{\sigma(D)} = q_{\sigma(C_1)}q_{\sigma(C_2)}\cdots q_{\sigma(C_r)} =
  q_{C_1}q_{C_2}\dots q_{C_r} = q_D.
  \]
  Now consider $i=(0_A,1_{V\setminus A})\in\cI$.  Then
  $\sigma(i)=(0_{\sigma(A)},1_{V\setminus\sigma(A)})$.  Using Proposition
  \ref{prop:satpara}, we obtain that
  \begin{align*}
    p_{\sigma(i)}&=\sum_{B:\sigma(A)\subseteq B} (-1)^{|B\setminus \sigma(A)|}
    q_B=\sum_{B:A\subseteq
    \sigma^{-1}(B)} (-1)^{|\sigma^{-1}(B)\setminus
    A|} q_{\sigma^{-1}(B)}= p_i.
  \end{align*}
\end{proof}

For a subset $C\subseteq V$, let $S(C) = \{ \sigma(C) \mid
\sigma\in\mathcal{S}\}$ be the orbit of $C$.
\begin{corollary}
  If the symmetry group $\mathcal{S}$ leaves the bi-directed graph $G$
  invariant, then the dimension of the marginal independence model with
  symmetry is
  \[
  \dim(\bP(G)\cap\bP(\mathcal{S}))=
  \sum_{\emptyset\not= C\in\cC(G)} 1/|S(C)|.
  \]
\end{corollary}
\begin{proof}
  By dividing through $|S(C)|$, every orbit of connected sets is counted once.
\end{proof}

\begin{corollary}
  \label{prop:symm-indep-1}
  If the symmetry group $\mathcal{S}$ leaves the bi-directed graph $G$
  invariant and the marginal independence model with symmetry is restricted
  to the interior of the probability simplex, then one obtains the curved
  exponential family $\bP_+(G)\cap\bP(\mathcal{S})$.
\end{corollary}
\begin{proof}
  The proof is analogous to the proof of Theorem 5.8.
\end{proof}

We define $\hat{n}_{\mathcal S} = n\cdot \hat{p}_{\mathcal S}$ to be the
fitted cell counts under the symmetry model ${\mathcal S}$, which are
simply the group averaged cell counts.  ML fitting of the model
$\bP(G)\cap\bP(\mathcal{S})$ may be performed by simply applying ICF for
fitting $\bP(G)$ to $\hat{n}_{\mathcal S}$ rather than the observed cell
counts. The rationale for this is as follows: let ${\mathcal L}_{G}(p;\{
n(i)\})$ indicate the likelihood for $\bP(G)$, evaluated with counts $\{
n(i)\}$.  If $p \in \bP(G)\cap\bP(\mathcal{S})$, then ${\mathcal
  L}_{G}(p;\{ n(i)\}) = {\mathcal L}_{G}(p;\{ \hat{n}_{\mathcal S}(i)\})$.
Further, ${\mathcal L}_{G}(p;\{ \hat{n}_{\mathcal S}(i)\}) = {\mathcal
  L}_{G}(\sigma(p);\{ \hat{n}_{\mathcal S}(i)\})$, for any $\sigma \in
{\mathcal S}$, where we define $\sigma(p(i)) = p(\sigma(i))$. Thus the
likelihood surface of the independence model $\bP(G)$ given the
group-averaged counts $\hat{n}_{\mathcal S}$ is invariant under
permutations $\sigma\in \mathcal{S}$ applied to probability vectors $p$. It
then follows that if $p^*$ is a local maximum of the likelihood function
${\mathcal L}_{G}(p;\{ \hat{n}_{\mathcal S}(i)\})$, then so is
$\sigma(p^*)$, for any $\sigma \in {\mathcal S}$. Further $p^*$ and
$\sigma(p^*)$ are in the same contour of the likelihood function.
Consequently if there is at most one local maximum of the likelihood
function ${\mathcal L}_{G}(p;\{ \hat{n}_{\mathcal S}(i)\})$ in any given
contour, then $p^* = \sigma(p^*)$ for all $\sigma \in {\mathcal S}$. Thus a
maximum found by ICF when applied to $\hat{n}_{\mathcal S}$, is in
$\bP(\mathcal{S})$, and is thus a maximum of the likelihood for the model
of symmetry and independence.

 \addtocounter{theorem}{-8}
\begin{example}\rm ({\it continued}).
Applying ICF to fit the model $A_1\ind D_2$ and $A_2\ind D_1$ for the twin data, using the fitted counts from the symmetry model $\mathcal{S}_{\hbox{\scriptsize twin}}$  resulted
in the fitted distribution shown in Table \ref{tab:twinsymmle}. 
The combined model has a deviance of $16.156$ on $2$ degrees of freedom,
taking the symmetry model given by $\mathcal{S}_{\hbox{\scriptsize twin}}$
as the alternative.  The corresponding p-value of $0.0003$ indicates a poor
fit and we may safely reject the generating hypothesis represented by the
graph in Figure \ref{fig:generating}(a).

  \end{example}
   \addtocounter{theorem}{-7}
 
   The approach taken here to combining symmetry and independence is
   analogous to that of
   \cite{andersson:madsen:symmetry:1998} in the Gaussian case. A more
   general approach would be to apply a symmetry group directly to the
   M\"{o}bius parameters, possibly with the restriction that orbits should
   only contain parameters corresponding to sets of a given cardinality;
   this would be more analogous to the work of
   \cite{hojsgaard:lauritzen:graphical:2006}.
  

\section{Related Work and Discussion}
\label{sec:conclu}

Several other authors have made use of the M\"{o}bius
decomposition or similar schemes.  \cite{lee:93} used this decomposition to
generate random binary vectors with fixed marginal distributions and
specified degrees of association.  \cite{ekholm:95,ekholm:00,ekholm:03}
used dependence ratios (see Remark~\ref{rem:dep-ratio}) to build
association and regression models for multivariate discrete responses.
Though Ekholm et al.~did not study marginal independence models {\em per se}, their work
on regression models offers one approach to building marginal independence
models for mixed continuous and discrete variables, which is an open
problem for future work.

Kauermann (\citeyear{kauermann:1997}) developed a parametrization for
marginal independence models using the multivariate logistic (m-logit)
transformation, which selects the highest order interaction term from every
margin.  However, the transformation from m-logit parameters to cell
probabilities cannot, in general, be computed in closed form.
Further, unlike classical log-linear parameters, the  valid m-logit
parameters may form a complicated subset of $\RRR^{2^V-1}$ and are not in
general variation independent.  The m-logit parameterization is a special
case of the  marginal log-linear framework of
\cite{bergsma:2002}.  In certain cases, such as for Figure
\ref{fig:bg}(a), there may exist a marginal log-linear parameterization for
a marginal independence model in which the parameters are variation
independent; see \cite{bergsma:2002},
\cite{lupparelli:marchetti:covariance:2005}.  However, there are models for
which this approach does not appear to lead to variation independent
parametrizations. Specifically, there does not appear to be such a
parametrization for the bi-directed chordless five cycle; see
\cite{bergsma:variation:2002} for related discussion.

As stated earlier, the problems inherent in expressing marginal
independence constraints in terms of a log-linear parametrization over a
larger set of variables are part of the general problem of `lack of upward
compatibility': specifically, a log-linear two-way interaction expresses a
property of the full joint distribution, and not of the relevant two-way
margin. A number of schemes have been proposed for dealing with this
problem, in addition to the m-logits mentioned above: see
\cite{ip:2003,streitberg:99,streitberg:90}. These provide alternative
parametrizations for the binary bi-directed models introduced here, which
may be computed from the fitted distribution, if desired.

\cite{cox:isi} and \cite{cox:94,coxwerm:book} take a different approach to
the problem of modelling independence structures similar to Gaussian
covariance models.  They focus on the quadratic binary exponential
distribution, also known as the Boltzmann machine \citep{hinton:83} or the
auto-logistic scheme \citep{bes:prague}.  In this distribution, the absence
of a given interaction term does not imply exact marginal independence, but
by approximating the marginal distributions via series expansions, it is
possible to gauge the size of any such dependence. As Cox notes, the extent
to which such marginal approximations are reasonable will depend on the
size of the relevant interaction terms.

\subsection*{Acknowledgments}

We would like to thank Steen Andersson, Sir David Cox, Steffen Lauritzen,
Fero Matu\v{s}, Alejandro Murua, Michael Perlman, James Robins, Tam\'{a}s
Rudas, Milan Stude{n}\'{y}, Jim Q.~Smith, Peter Spirtes, Bernd Sturmfels
and Nanny Wermuth for helpful conversations.  Kenneth Kendler provided the
twin data.  Yen-Sheng Chiang and Richard Callahan suggested the analysis of
questions relating to Trust.  We are particularly grateful to Anders Ekholm
for his comments on earlier versions of this paper.  This research was
supported by the U.S.\ National Science Foundation (DMS-9972008,
DMS-0505612, DMS-0505865), the U.S.\ National Institutes for Health
(R01-HG2362-3), the William and Flora Hewlett Foundation and the Center for
Advanced Studies in the Behavioral Sciences at Stanford University where
Thomas Richardson was a Fellow from 2003-2004.

\appendix
\setcounter{section}{1} 

\section*{Appendix: Likelihood Equations and Hessian calculations}

If $\G$ is a bi-directed graph with vertex set $V$, then for an
arbitrary subset $A\subseteq V$, let 
\[
[A]_{\G} = \{ C \mid C \text{ is a
maximal connected component of } \G_A\}.
\]
Note that $[A]_\G$ forms a partition $A=\bigcup_{C\in [A]_\G} C$.  For
disconnected sets $D\subseteq V$ this partition is the one used in Theorem
\ref{thm:moebconstraints}.  Since for a connected set $C\subseteq V$ the
family $[C]_\G$ only comprises one set, namely $C$ itself, we have that
under a joint distribution in the model $\bP(\G)$,
\[
q_A = \prod_{C \in [A]_{\G}} q_C,\quad A\subseteq V.
\]
Hence for any set $A$, there is a unique expansion of the joint cell
probability $p_A^V$ in terms of the parameters $q_C$ for connected sets
$C$ in $\G$,
\[
p_A^V = \sum_{B : A\subseteq B} (-1)^{|B\setminus A|} \prod_{C: C\in
  [B]_{\G}} q_C,
\]
recall that $p^V_A = P(X_A=0,X_{V\setminus A}=1)$. We call this last expression the {\it expansion} for $p_A^V$ (under graph $\G$).

\begin{lemma}\label{lem:expansion}
If $C$ is a connected set in the graph $\G$, then the
 parameter $q_C$ appears in the expansion for $p_A^V$ if and only if
$\spo(C)\cap (A\setminus C) = \emptyset$.
\end{lemma}
\begin{proof}
If $\spo(C)\cap (A\setminus C) = \emptyset$
then $C \cup (A\setminus C)$ forms a disconnected superset of $A$
in which $C$ is a maximal connected component. Hence 
$C\in [C\cup (A\setminus C)]_{\G}$.
If $\spo(C)\cap (A\setminus C) \neq \emptyset$ then there is a vertex $a
\in A\setminus C$ such that $a\in \spo (C)$.  Hence, in any set $B$
containing $A$ and $C$, there is a maximal connected set $\bar C \supseteq
C\cup\{a\}$. Hence $C\notin [B]_\G$ for any $B\supseteq A$.
\end{proof}

In words, Lemma \ref{lem:expansion} states that parameter $q_C$ appears in
the expansion for $p_A^V$ if and only if every vertex in $A$ that is
adjacent to $C$ is already in $C$.  Consequently, $(\partial/\partial q_C)
p_A^V = 0$ for any connected set $C$ in $\G$ that satisfies $\spo (C) \cap
(A\setminus C) \neq \emptyset$.

\begin{lemma}\label{derivative}
  If $\spo (C) \cap (A\setminus C) = \emptyset$ then
\[
{{\partial p_A^V}\over {\partial q_C}} = (-1)^{|C\setminus A|}
p_{A\setminus C}^{V\setminus \spo(C)}.
\]
\end{lemma}

\begin{proof}  The claim holds since $\spo(C)\cap (A\setminus C) = \emptyset$ iff
  $(A\setminus C) \subseteq (V\setminus \spo(C))$, and
\begin{eqnarray*}
{{\partial p_A^V}\over {\partial q_C}} &=&
\sum_{\scriptstyle
B : (A\setminus C) \subseteq B\atop\scriptstyle
\mathrm{ and } B\subseteq V\setminus \spo(C)
}
(-1)^{|(B \dot{\cup} C)\setminus A|}
\prod_{C^*:C^*\in [B]_{\G}} q_{C^*}\\
&=& (-1)^{|C\setminus A|} \sum_{\scriptstyle
B : (A\setminus C) \subseteq B\atop\scriptstyle
\mathrm{ and } B\subseteq V\setminus \spo(C)}
(-1)^{|B \setminus A|}
\prod_{C^*:C^*\in [B]_{\G}} q_{C^*}\\
&=& (-1)^{|C\setminus A|} p_{A\setminus C}^{V\setminus \spo(C)}.\\
\end{eqnarray*}
\vskip-40pt
\end{proof}

\begin{corollary}\label{cor:likelihoodeqns}
  The system of likelihood equations associated with the model $\bf{P}(\G)$
  contains an equation
\[
{{\partial \ell}\over {\partial q_C}} =
\sum_{A: \spo (C) \cap (A\setminus C) = \emptyset}
(-1)^{|C\setminus A|} {n_A^V \over p_A^V} p_{A\setminus C}^{V\setminus \spo(C)} = 0
\]
for every (non-empty) connected set $C$ in $\G$.
\end{corollary}

The likelihood equations can also be expressed in terms of expectations
with respect to conditional empirical measures (provided these exist):
\[
\mathbb{E}_{X_{V\setminus (\spo(C) \setminus C)} \mid X_{\spo(C)
    \setminus C}=1} \left[ (-1)^{\sum_{i\in C} X_i}P\left( X_{\spo (C)}
    \mid X_{V\setminus \spo(C)}\right)^{-1}\right] = 0
\]
where $\mathbb{E}_{X_{V\setminus (\spo(C) \setminus C)} \mid
  X_{\spo(C) \setminus C}=1}$ is expectation w.r.t.~the measure on
$\cI_{V\setminus (\spo(C)\setminus C)}$ given by (normalizing)
the empirical frequencies in the sub-table in which $X_{\spo(C)
  \setminus C}=1$.

\begin{lemma}
  Let $C$ and $\bar C$ be connected sets in $\G$.
  \begin{itemize}
  \item[(i)] If $\left(\spo(C)\cap (A\setminus C)\right)\cup
    \left(\spo(\bar C)\cap (A\setminus \bar C)\right)\neq\emptyset$, then
    the second derivative
    \[
    {\partial \over{\partial q_C}}{\partial \over{\partial q_{\bar C}}} \log
    p_A^V =0.
    \]
  \item[(ii)] If $\left(\spo(C)\cap (A\setminus C)\right)\cup
    \left(\spo(\bar C)\cap (A\setminus \bar C)\right)=\emptyset$ and $\bar
    C \cap \spo(C) \neq \emptyset$, then
    \[
    {\partial \over{\partial q_C}}{\partial \over{\partial q_{\bar C}}} \log
    p_A^V 
    = - (-1)^{|C\setminus A|}(-1)^{|\bar C\setminus A|}
    p_{A\setminus C}^{V\setminus\spo(C)}p_{A\setminus \bar
    C}^{V\setminus\spo(\bar C)}\cdot {1\over {(p_A^V)^2}}. 
    \]
  \item[(iii)] If  $\left(\spo(C)\cap (A\setminus C)\right)\cup
    \left(\spo(\bar C)\cap (A\setminus \bar C)\right)=\emptyset$ and
    $\bar C \cap \spo(C) = \emptyset$, then
    \begin{multline*}
      {\partial \over{\partial q_C}}{\partial \over{\partial q_{\bar C}}}
      \log p_A^V =(-1)^{|(C\cup \bar C) \setminus A)|}
      p_{A\setminus (C\cup \bar C)}^{V\setminus \spo (C\cup \bar C)}\cdot
      {1\over{p_A^V}} 
\\
      - (-1)^{|C\setminus
        A|}(-1)^{|\bar C\setminus A|} p_{A\setminus C}^{V\setminus 
        \spo(C)}p_{A\setminus \bar C}^{V\setminus 
        \spo(\bar C)}\cdot {1\over {(p_A^V)^2}}.
    \end{multline*}
  \end{itemize}
\end{lemma}
\begin{proof} This follows from Lemma \ref{derivative}.
  The first term on the RHS of the equation in (iii) occurs if the derivative of $(\partial/\partial
  \bar C) p_{A\setminus C}^{V\setminus \spo (C)}$ is non-zero,
  which requires $\bar C\subseteq V\setminus \spo (C)$ and $\spo
  (\bar C) \cap \left( A\setminus (C\cup \bar C)\right) = \emptyset$.  The
  second condition is implied by $\spo (\bar C) \cap (A\setminus \bar C)
  =\emptyset$.  The first is equivalent to $\bar C \cap \spo(C) =
  \emptyset$.
\end{proof}
In words, the condition that $\bar C \cap  \spo(C)=\emptyset$
requires that $C$ and $\bar C$ are disjoint and there
is no vertex in $C$ adjacent to a vertex in $\bar C$.
Note that $\spo (C) \cap \bar C = \emptyset$ if and only if 
$\spo (\bar C) \cap C =\emptyset$, hence the conditions in (ii) and (iii)
are symmetric in $C$ and $\bar C$ as required.

The full Hessian may be obtained by summing the expression given in the
last Lemma over all sets $A\subseteq V$.

\small
\bibliographystyle{Chicago}
\bibliography{bmmi}
\end{document}